\documentclass[preprint,12pt]{elsarticle}

\usepackage{graphicx}

\usepackage{amssymb}
\usepackage[fleqn]{amsmath}
\usepackage{amsthm}

\newtheorem{thm}{Theorem}

\theoremstyle{example}

\theoremstyle{definition}

\theoremstyle{remark}

\journal{}

\begin{document}

\begin{frontmatter}



\title{Limiting case Hardy inequalities on the sphere}

\author{Ahmed A. Abdelhakim}
\address{Mathematics Department, Faculty of Science, Assiut University, Assiut 71516 - Egypt\\
Email: ahmed.abdelhakim@aun.edu.eg}
\begin{abstract}
We give sharp limiting case Hardy inequalities on the sphere $\mathbb{S}^{2}$ and show that their optimal constants are unattainable by any
$f\in H^{1}\left(\mathbb{S}^{2}\right)\setminus\{0\}$. The singularity of the problem is related to the geodesic distance from a point on the sphere.
\end{abstract}

\end{frontmatter}
\section{Introduction}
The classical Hardy inequality
\begin{equation}\label{ie1}
\int_{\mathbb{R}^n} |\nabla u|^{2} dx
\geq \frac{(n-2)^{2}}{4}\int_{\mathbb{R}^n} \frac{u^{2}}{|x|^2} dx
\end{equation}
is valid in dimensions $n\geq 3$
for all functions $u\in H^{1}\left(\mathbb{R}^{n}\right)$ (\cite{Balinsky}).
It obviously fails on $\mathbb{R}^{2}$ as the right hand side of (\ref{ie1}) no longer makes sense.
In order to obtain a version of (\ref{ie1})
in the critical case $n=2$ on bounded domains, a logarithmic weight can be introduced to tame the singularity. In \cite{Ioku1,Ioku2,Machihara,
Ruzhansky,Sano,Takahashi1,Takahashi2},
for instance, inequalities of the type
\begin{equation*}
\int_{B} |\nabla u|^{n} dx
\geq C_{n}(\Omega)\int_{B} \frac{|u|^{n}}{|x|^n
\left(\log{\frac{1}{|x|}}\right)^{n}} dx
\end{equation*}
were analysed for $u \in W_{0}^{1,n}(B)$ where $B$ is the unit ball in $\mathbb{R}^{n}$.\\
\indent
Let $n\geq 3$ and $\,\mathbb{S}^{n}$ be the unit sphere equipped with its Lebesgue surface measure
$\sigma_{n}$ in $\mathbb{R}^{n+1}$. Denote by $d(.,p):\mathbb{S}^{n}\rightarrow [0,\pi]$ the geodesic distance from $p\in \mathbb{S}^{n}$, and by $\nabla_{\mathbb{S}^{2}}$ the gradient on $\mathbb{S}^{n}$. Recently, Xiao \cite{Xiao} proved that
if $f\in C^{\infty}\left(\mathbb{S}^{2}\right)$ then
 \begin{equation}\label{xiao}
{\bar{c}_{n}}\int_{\mathbb{S}^{n}} f^{2} d \sigma_{n}+
\int_{\mathbb{S}^{n}} |\nabla_{\mathbb{S}^{2}} f|^{2}
d\sigma_{n} \geq c_{n}^{2}
\int_{\mathbb{S}^{n}}
\left(\frac{f^{2}}{d\left(x,p\right)^{2}}+
\frac{f^{2}}{\left(\pi-d(x,p)\right)^{2}} \right) d \sigma_{n}
\end{equation}
with $ \bar{c}_{n}
=\left(\frac{2}{3}+\frac{1}{\pi^2}\right)
c_{n}^{2}+c_{n},\,$
$c_{n}=\frac{n-2}{2}.$
It was also shown in \cite{Xiao} that the constant
$c_{n}$ in (\ref{xiao}) is sharp in the sense that
\begin{equation*}
c_{n}^{2}=
\inf_{f\in C^{\infty}\left(\mathbb{S}^{n}\right)\setminus\{0\}}
\frac{D_{n}(f)}
{\int_{\mathbb{S}^{n}}
\frac{f^{2}}{d\left(x,p\right)^{2}}d \sigma_{n}} \,=\,
\inf_{f\in C^{\infty}\left(\mathbb{S}^{n}\right)\setminus\{0\}}
\frac{D_{n}(f)}
{\int_{\mathbb{S}^{n}}
\frac{f^{2}}{\left(\pi-d(x,p)\right)^{2}} d \sigma_{n}}
\end{equation*}
where
\begin{equation*}
D_{n}(f):={c_{n}}\int_{\mathbb{S}^{n}} f^{2} d \sigma_{n}+
\int_{\mathbb{S}^{n}} |\nabla_{\mathbb{S}^{2}} f|^{2}
d\sigma_{n},\quad f\in C^{\infty}\left(\mathbb{S}^{n}\right).
\end{equation*}
\indent We prove $L^{2}$ Hardy inequalities with  optimal constants on the sphere $\mathbb{S}^{2}$ in $\mathbb{R}^{3}$. This is a critical exponent case
as the integral $\int_{\mathbb{S}^{2}}{\theta^{-1+\lambda}}
{d\sigma_{2}}$, where $\theta$ is the polar angle, diverges for $\lambda \leq -1$.
We also argue the lack of maximizers for
our inequalities. Our approach denies the possibility
of an equality in Xiao's inequality (\ref{xiao}) as well.
\section{Preliminaries}
A point on the sphere $\mathbb{S}^{2}$
will have the standard spherical coordinate parametrization $\left(\sin{\theta}\cos{\varphi},\sin{\theta}\sin{\varphi},
\cos{\theta}\right)$ where  $\theta \in [0,\pi]$ refers to the polar angle and $\varphi\in [0,2\pi[$ is the azimuthal angle. Then the surface measure induced by the Lebesgue measure on $\mathbb{R}^3$ is $d\sigma_{2}=\sin{\theta} d\theta d\varphi$, the gradient and
the Laplace-Beltrami operator, respectively, are given by
\begin{equation*}
\nabla_{\mathbb{S}^{2}}=
\hat{\theta}\,\frac{\partial}{\partial \theta}+\hat{\varphi}\,\frac{1}{\sin{\theta}}\frac{\partial}{\partial \varphi},\;
\Delta_{\mathbb{S}^{2}}=\frac{1}{\sin{\theta}}
\frac{\partial}{\partial \theta}\left(\sin{\theta}\,\frac{\partial}{\partial \theta}\right)+\frac{1}{\sin^{2}{\theta}}
\frac{\partial^{2}}{\partial \varphi^{2}}.
\end{equation*}
The Sobolev space $H^{1}\left(\mathbb{S}^{2}\right)$
is the completion of $C^{\infty}\left(\mathbb{S}^{2}\right)$ in the norm
\begin{equation*}
\parallel f\parallel_{H^{1}\left(\mathbb{S}^{2}\right)}:=
\left(\parallel f\parallel^{2}_{L^{2}\left(\mathbb{S}^{2}\right)}+
\parallel \nabla f \parallel^{2}_{L^{2}\left(\mathbb{S}^{2}\right)}
\right)^{\frac{1}{2}}.
\end{equation*}
\indent
In order to find the geodesic distance $d(x,p)$ from a point $x\in \mathbb{S}^{2}$ to a given a point $p\in \mathbb{S}^{2}$, we rotate the axes,
if necessary, to put $p$ on the zenith direction
then place the great circle passing through $p$ and $x$ in the azimuth reference direction so that we have $d(x,p)=\theta$.\\
\indent For simplicity, we henceforth denote
$d\sigma_{2}$, $\nabla_{\mathbb{S}^{2}}$ and $\Delta_{\mathbb{S}^{2}}$ by $d\sigma$, $\nabla$ and $\Delta$, respectively.
\section{Main results}
Let $\phi:[0,\pi]\rightarrow [1,\infty[$ be defined by $\phi(t):=\log\left(\pi e/t\right),\,$
$\psi:[0,\pi]\rightarrow [1+\log{\pi},\infty[$
be such that $\psi(t):=\phi\left(\sin{t}\right),$
and $\,\rho_{\phi}(t):=t\phi(t)$. Let $A>0$.
Denote by $S,$ $\,T_{A},$ and $\,Q\left(.;\phi\right)$
the positive nonlinear functionals on
$H^{1}\left(\mathbb{S}^{2}\right)$
given by
\begin{align*}
S(f):=&\int_{\mathbb{S}^{2}}
|\hat{\theta}.\nabla f|^{2} d\sigma+\frac{1}{2\pi^{2}}\int_{\mathbb{S}^{2}}
f^{2}\,d\sigma,\\
T_{A}(f):=&\int_{\mathbb{S}^{2}} |\nabla  \, f|^{2}
d\sigma_{2}+\frac{A}{4}\int_{\mathbb{S}^{2}} f^{2} d \sigma_{2},\\
Q\left(f;\phi\right):=&
\frac{1}{4}
\int_{\mathbb{S}^{2}}
\bigg(
\frac{f^{2}}{{{\rho}_{\phi}}^{2}\left(d\left(x,p\right)\right)}+
\frac{f^{2}}{{{\rho}_{\phi}}^{2}\left(\pi-d\left(x,p\right)\right)}\bigg)
d \sigma_{2}.
\end{align*}
\begin{thm}\label{thm1}
Assume that $f\in H^{1}\left(\mathbb{S}^{2}\right).$ Then there exists constants $A,\,B>0$, independent of $f$, such that
\begin{eqnarray}
\label{main11}  Q\left(f;\phi\right)&\leq& T_{A}(f),\\
\label{main12}  Q\left(f;\psi\right)&\leq& T_{B}(f).
\end{eqnarray}
\end{thm}
Both inequalities (\ref{main11})
and (\ref{main12}) are optimal, but an equality
is impossible in either one:
\begin{thm}\label{thm2}
\begin{gather}
\label{main21}
\sup_{f\in H^{1}\left(\mathbb{S}^{2}\right)\setminus\{0\}}
\frac{Q\left(f;\phi\right)}
{T_{A}(f)} = 1, \\
\label{main22}
\sup_{f\in H^{1}\left(\mathbb{S}^{2}\right)\setminus\{0\}}
\frac{Q\left(f;\psi\right)}
{T_{B}(f)} = 1.
\end{gather}
\end{thm}
\begin{thm}\label{thm3}
There does not exist $f\in H^{1}\left(\mathbb{S}^{2}\right)\setminus\{0\}$ such that $\,Q\left(f;\phi\right)= T_{A}(f),$
or $\,Q\left(f;\psi\right)= T_{B}(f)$.
\end{thm}
A variant of the abovementioned results follows via a different approach:
\begin{thm}\label{thm4}
Let
$f\in H^{1}\left(\mathbb{S}^{2}\right).$ Then
\begin{align}
\label{w01}
&\frac{1}{4}\int_{\mathbb{S}^{2}} \frac{f^{2}}{{{\rho}_{\phi}}^{2}\left(d\left(x,p\right)\right)}
\, d\sigma
\leq S(f)+\frac{1}{2\pi}\int_{\mathbb{S}^{2}}
\frac{f^{2}}{\pi-d\left(x,p\right)}\, d\sigma,\\
\label{w02}
&\frac{1}{4}
\int_{\mathbb{S}^{2}} \frac{f^{2}}{{{\rho}_{\phi}}^{2}\left(\pi-d\left(x,p\right)\right)}
\, d\sigma \leq S(f)+\frac{1}{2\pi}\int_{\mathbb{S}^{2}}
\frac{f^{2}}{d\left(x,p\right)}\, d\sigma.
\end{align}
Moreover
\begin{equation}\label{w03}
\sup_{f\in H^{1}\left(\mathbb{S}^{2}\right)\setminus\{0\}}
\frac{\frac{1}{4}\int_{\mathbb{S}^{2}} \frac{f^{2}}{{{\rho}_{\phi}}^{2}\left(d\left(x,p\right)\right)}
\, d\sigma}
{S(f)+\frac{1}{2\pi}\int_{\mathbb{S}^{2}}
\frac{f^{2}}{\pi-d\left(x,p\right)}\, d\sigma} =
\sup_{f\in H^{1}\left(\mathbb{S}^{2}\right)\setminus\{0\}}
\frac{\frac{1}{4}\int_{\mathbb{S}^{2}}\frac{f^{2}}{{{\rho}_{\phi}}^{2}\left(\pi-d\left(x,p\right)\right)}
\, d\sigma}
{S(f)+\frac{1}{2\pi}\int_{\mathbb{S}^{2}}
\frac{f^{2}}{d\left(x,p\right)}\, d\sigma} =
1,
\end{equation}
and the suprema in (\ref{w03}) are not attained
in $H^{1}\left(\mathbb{S}^{2}\right)\setminus\{0\}$.
\end{thm}
\section{Proof of Theorem 1}
\begin{proof}
Let $f\in C^{\infty}\left(\mathbb{S}^{2}\right)$.
Notice that $\psi>1$ and write $f(\theta,\varphi)=\sqrt{\psi(\theta)} g(\theta,\varphi)$. We have
\begin{eqnarray}
\nonumber
|\nabla  \,  f|^{2}&=& |\psi^{\frac{1}{2}} \nabla  \,  g+
 g \nabla  \,  \psi^{\frac{1}{2}}|^{2}\\
\nonumber &=& \psi |\nabla  \,  g|^{2} +\langle\psi^{\frac{1}{2}} \nabla  \,  g, g\psi^{-\frac{1}{2}} \nabla  \,\,\psi\rangle+ |\frac{1}{2}
{\psi}^{-\frac{1}{2}}\nabla  \,\psi|^{2} g^{2}\\
\label{p1}&=& \psi |\nabla  \,  g|^{2}
 +\frac{1}{2}\langle
\nabla  \,\psi,\,
\nabla  \,  g^2\rangle+\frac{1}{4}
\frac{1}{\psi}
|\nabla  \,\psi|^{2} g^{2}.
\end{eqnarray}
Integrating both sides of (\ref{p1}) over $\mathbb{S}^{2}$ we get
\begin{eqnarray}
\nonumber
\int_{\mathbb{S}^{2}} |\nabla  \,  f|^{2} d\sigma
&=& \int_{\mathbb{S}^{2}} \left(  \psi |\nabla  \,  g|^{2}
 +\frac{1}{2}\langle
\nabla  \,\psi,\,
\nabla  \,  g^2\rangle+\frac{1}{4}
\frac{1}{\psi}
|\nabla  \,\psi|^{2} g^{2}
 \right) d\sigma\\
\label{notattain} &\geq& \frac{1}{4}\int_{\mathbb{S}^{2}}
\frac{1}{\psi}
|\nabla  \,\psi|^{2} g^{2}
  d\sigma
  +\frac{1}{2}\int_{\mathbb{S}^{2}} \langle
\nabla  \,\psi,\,
\nabla  \,  g^2\rangle
 d\sigma\\
 \label{p2} &=& \frac{1}{4}\int_{\mathbb{S}^{2}}
\frac{1}{\psi}
  |\psi^{\prime}|^{2} g^{2}
  d\sigma
  -\frac{1}{2}\int_{\mathbb{S}^{2}}
  g^{2}  \Delta \psi
 d\sigma
\end{eqnarray}
by partial integration over the closed manifold $\mathbb{S}^{2}$. Calculating, we find
\begin{align}\label{delt}
\Delta \psi=
\frac{1}{\sin{ \theta}}\frac{\partial}{\partial \theta}
\left( \sin{\theta}\frac{\partial }{\partial \theta}\psi \right)=
1.
\end{align}
Returning $g$ to $f/\sqrt{\psi}$ and substituting for $\Delta \psi$ from (\ref{delt}) into (\ref{p2}), we obtain
\begin{eqnarray}
\label{p3}
\int_{\mathbb{S}^{2}} |\nabla  \,  f|^{2} d\sigma &\geq&  \frac{1}{4}\int_{\mathbb{S}^{2}}
 \frac{f^{2}}{\psi^{2}}
 \frac{\cos^{2}{\theta}}{\sin^{2}{\theta}}
d\sigma- \frac{1}{2}
\int_{\mathbb{S}^{2}}
 \frac{f^{2}}{\psi}
 d\sigma.
\end{eqnarray}
Adding the finite integral $\,\displaystyle
\frac{1}{4}\int_{\mathbb{S}^{2}}
\left(\frac{1}{\theta^{2}\phi^{2}\left(\theta\right)}+
\frac{1}{\left(\pi-\theta\right)^{2}\phi^{2}
\left(\pi-\theta\right)}\right)f^{2}
d\sigma\,$ to both sides of (\ref{p3})
transforms it into the inequality
\begin{align}
&\nonumber \frac{1}{4}\int_{\mathbb{S}^{2}}
\left(\frac{1}{\theta^{2}\phi^{2}\left(\theta\right)}+
\frac{1}{\left(\pi-\theta\right)^{2}\phi^{2}
\left(\pi-\theta\right)}\right)f^{2}
d\sigma\\
\label{p4}&\qquad\qquad\qquad\leq\int_{\mathbb{S}^{2}} |\nabla  \,  f|^{2} d\sigma+\frac{1}{4}\int_{\mathbb{S}^{2}} F(\theta)\,  f^{2} d\sigma,
\end{align}
where
\begin{equation*}
F(t):=
\frac{1}{t^{2}\phi^{2}\left(t\right)}+
\frac{1}{\left(\pi-t\right)^{2}\phi^{2}
\left(\pi-t\right)}-
\frac{\cos^{2}{t}}{\sin^{2}{t}}
\frac{1}{\phi^{2}\left(\sin{t}\right)}
+\frac{2}{\phi\left(\sin{t}\right)}.
\end{equation*}
Obviously, $F$ is continuous on $]0,\pi[$
and, as expected from the facts that
$\phi(t)\rightarrow +\infty$ when $t\rightarrow 0^{+},\,$ $\,\sin{t} = t+o(t)$ as $t\rightarrow 0,\,$
it turns out
\begin{equation*}
\lim_{t\rightarrow 0^{+}} F(t)=
\lim_{t\rightarrow \pi^{-}} F(t)=\frac{1}{\pi^{2}}.
\end{equation*}
Hence, $F$ can be extended to a uniformly continuous, consequently a bounded, function on $[0,\pi]$. Noting this in (\ref{p4}) implies (\ref{main11}). Direct computation also shows
\begin{equation*}
A=\sup_{[0,\pi]}|F|=F(\frac{\pi}{2})=
\frac{2}{1+\log{\pi}}+\frac{8}{\left(1+\log{2}\right)^{2}}
\frac{1}{\pi^{2}}.
\end{equation*}
To prove (\ref{main12}), we add
to both sides of (\ref{p3})
the well-defined  integral \\\\ $\,\displaystyle
\frac{1}{4}\int_{\mathbb{S}^{2}}
\left(\frac{1}{\theta^{2}}+
\frac{1}{\left(\pi-\theta\right)^{2}}\right)
\frac{f^{2}}{\psi^{2}\left(\theta\right)}
d\sigma.$ We then obtain the following analogue of
(\ref{p4}):
 \begin{align}
&\nonumber \frac{1}{4}\int_{\mathbb{S}^{2}}
\left(\frac{1}{\theta^{2}}+
\frac{1}{\left(\pi-\theta\right)^{2}}\right)
\frac{f^{2}}{\psi^{2}\left(\theta\right)}
d\sigma\\
\label{p5}&\qquad\qquad\qquad\leq\int_{\mathbb{S}^{2}} |\nabla  \,  f|^{2} d\sigma+\frac{1}{4}\int_{\mathbb{S}^{2}} G(\theta)\,  f^{2} d\sigma,
\end{align}
where
\begin{eqnarray}
\nonumber
G(t)&:=&
\frac{M(t)}{\psi^{2}(t)}
+\frac{2}{\psi(t)},\\
\label{mm1}
M(t)&:=&\frac{1}{t^{2}}+
\frac{1}{(\pi-t)^{2}}
-\frac{\cos^{2}{t}}{\sin^{2}{t}}.
\end{eqnarray}
Once the boundedness
of $G$ is ensured, we see that (\ref{p5})
yields the inequality (\ref{main12}).
Evidently, $G$ has the same features as $F$.
Since \begin{equation}\label{mm2}
\lim_{\theta\rightarrow 0} M(\theta)
\,=\,\lim_{\theta\rightarrow \pi} M(\theta)
 \,=\,\frac{2}{3}+\frac{1}{\pi^2},\;
\lim_{\theta\rightarrow 0^{+}} \psi(t)
\,=\,\lim_{\theta\rightarrow \pi^{-}} \psi(t)
\,=\,+\infty
\end{equation}
then $\,M\in C[0,\pi],$ and
$\;\lim_{t\rightarrow 0^{+}} G(t)=
\lim_{t\rightarrow \pi^{-}} G(t)=0,\,$
which makes $G$ bounded on $[0,\pi]$. Moreover
\begin{equation*}
B=\sup_{[0,\pi]}|G|=G(\frac{\pi}{2})=
\frac{2}{1+\log{\pi}}+\frac{8}{\left(1+\log{\pi}\right)^{2}}
\frac{1}{\pi^{2}}.
\end{equation*}
\end{proof}
\section{Proof of Theorem \ref{thm2}}\label{proofthm2}
\begin{proof}
First, we would like to define the weak
laplace-Beltrami gradient
of a function $f\in L^{1}\left(\mathbb{S}^{2}\right)$.
Suppose $f\in C^{\infty}\left(\mathbb{S}^{2}\right)$
and $\,v(\theta,\varphi)=v_{\theta}(\theta,\varphi)
\hat{\theta}+v_{\varphi}(\theta,\varphi)\hat{\varphi}\,$
with $v_{\theta},v_{\varphi}\in  C^{\infty}\left(\mathbb{S}^{2}\right)$.
Then
\begin{align*}
&\int_{\mathbb{S}^{2}} \frac{\partial f}{\partial \theta}
\,v_{\theta} d\sigma =
\int_{\mathbb{S}^{2}}
\nabla   f \cdot \hat{\theta}
\,v_{\theta}d\sigma=-\int_{\mathbb{S}^{2}}
 f\, \nabla  \cdot (\,v_{\theta}\hat{\theta})
 d\sigma,
 \\
&\int_{\mathbb{S}^{2}}\frac{1}{\sin{\theta}} \frac{\partial f}{\partial \varphi}
\,v_{\varphi} d\sigma =
\int_{\mathbb{S}^{2}}
\nabla   f \cdot \hat{\varphi}
\,v_{\varphi}d\sigma=-\int_{\mathbb{S}^{2}}
 f\, \nabla  \cdot (\,v_{\varphi}\hat{\varphi})
 d\sigma.
\end{align*}
Adding these identities we get
\begin{equation}\label{mot1}
\int_{\mathbb{S}^{2}}\nabla f \cdot V \, d\sigma
=-\int_{\mathbb{S}^{2}} f\, \nabla\cdot V \, d\sigma
\end{equation}
for any vector field
$V\in C^{\infty}\left(\mathbb{S}^{2}\rightarrow
T\left(\mathbb{S}^{2}\right)\right)$ where $T\left(\mathbb{S}^{2}\right)$
is the tangent bundle of the smooth manifold $\mathbb{S}^{2}$. Motivated by (\ref{mot1}),  $f$ is weakly differentiable if there is
a vector field $\vartheta_{f}\in L^{1}\left(\mathbb{S}^{2}\rightarrow T\left(\mathbb{S}^{2}\right)\right)$ such that
\begin{equation}\label{mot2}
\int_{\mathbb{S}^{2}} \vartheta_{f} \cdot V \, d\sigma
=-\int_{\mathbb{S}^{2}} f\, \nabla\cdot V \, d\sigma,
\quad \forall\, V\in C^{\infty}\left(\mathbb{S}^{2}\rightarrow
T\left(\mathbb{S}^{2}\right)\right).
\end{equation}
This, unique up to a set of zero measure,  vector field $\vartheta_{f}$ is the weak surface gradient of $f$.
According to (\cite{Eichhorn}, Proposition 3.2., page 15)
\begin{equation*}
H^{1}\left(\mathbb{S}^{2}\right)=
W^{1,2}\left(\mathbb{S}^{2}\right):=
\left\{f\in L^{2}(\mathbb{S}^{2}):
|\vartheta_{f}|\in L^{2}\left(\mathbb{S}^{2}\right)\right\}.
\end{equation*}
We start with (\ref{main21}). By Theorem \ref{thm1}, it suffices to prove
the existence of a sequence $\left\{f_{n}\right\}_{n\geq 1}$ in $ H^{1}\left(\mathbb{S}^{2}\right)$
such that
\begin{equation}\label{limq}
\lim_{n\rightarrow \infty}\frac{Q\left(f_{n};\phi\right)}
{T_{A}(f_{n})}=1.
\end{equation}
Consider the functions
\begin{equation}\label{fn}
f_{n}(\theta,\varphi):=
{\phi(\theta)}^{\frac{1}{2}-\frac{1}{n}}.
\end{equation}
The functions $f_{n}$ are independent of $\varphi$,
hence
\begin{equation}\label{plug}
\frac{Q\left(f_{n};\phi\right)}
{T_{A}(f_{n})}=
\frac{
\int_{0}^{\pi}
\frac{f_{n}^{2}\,\sin{\theta}}{\theta^{2}\,\phi^{2}(\theta) }\,d\theta+
\int_{0}^{\pi}
\frac{f_{n}^{2}\,\sin{\theta}}
{\left(\pi-\theta\right)^{2}
\,\phi^{2}\left(\pi-\theta\right)
}\,d\theta}
{4\int_{0}^{\pi} \left(\frac{\partial f_{n}}{\partial \theta} \right)^{2}\,\sin{\theta}
d\theta+A \int_{0}^{\pi} f_{n}^{2} \,\sin{\theta}
d\theta}
\end{equation}
where the derivative $\partial f_{n}/\partial \theta$
is understood in the week sense discussed above.
Since $\phi \in L^{1}_{\text{loc}}
\left(\mathbb{R}\right)$ and $\phi \geq 1$
on $[0,\pi]$, then
\begin{eqnarray}
\label{calc1}
\int_{0}^{\pi} f^{2}_{n} \,\sin{\theta}
d \theta\,=\,
\int_{0}^{\pi} {\phi(\theta)}^{1-\frac{2}{n}} \,\sin{\theta}
d\theta
\leq \int_{0}^{\pi} {\phi(\theta)} \,d\theta
\,\approx\, 1.
\end{eqnarray}
Thus $f_{n} \in L^{2}
\left(\mathbb{S}^{2}\right)$ for all $n\geq 1$.
Notice also that $f_{n}$ is smooth on $[0,\pi]\setminus \{0\}$ and its weak derivative
\begin{equation}\label{fnprime}
\frac{\partial f_{n}}{\partial \theta} =\frac{\frac{1}{n}-\frac{1}{2}}
{\theta\,{\phi^{\frac{1}{2}+\frac{1}{n}}}}.
\end{equation}
Therefore
\begin{eqnarray}
\nonumber
\int_{0}^{\pi} \left(\frac{\partial f_{n}}{\partial \theta}\right)^{2} \,\sin{\theta}
d \theta\,=\,
\frac{a_{n}}{4}
\int_{0}^{\pi}
\frac{1}
{\theta\,{\phi^{1+\frac{2}{n}}}}
 \,\frac{\sin{\theta}}{\theta}
d \theta,\quad a_{n}:=\left(1-\frac{2}{n}\right)^{2}.
\end{eqnarray}
And since $\displaystyle \,\int_{0}^{\pi}
\frac{d \theta}{\theta\,{\phi^{1+\frac{2}{n}}}}
 =\frac{n}{2},\,$
$\,{\sin{\theta}}\leq {\theta},\,$
then $\,{\partial f_{n}}/{\partial \theta} \in L^{2}\left(\mathbb{S}^{2}\right)$
for all $n\geq 1$. Substituting for $f_{n}$ from (\ref{fn}) and for ${\partial f_{n}}/{\partial \theta}$ from (\ref{fnprime}) into
(\ref{plug}) implies
\begin{equation}\label{lim}
\frac{Q\left(f_{n};\phi\right)}
{T_{A}(f_{n})}=
\frac{\alpha_{n}+\beta_{n}}{a_{n}\alpha_{n}+\gamma_{n}}=
\frac{1}{a_{n}}\left(
1+\frac{\beta_{n}-
\gamma_{n}/a_{n}}{\alpha_{n}+\gamma_{n}/a_{n}}\right)
\end{equation}
where
\begin{eqnarray*}
  \alpha_{n} &:=&\int_{0}^{\pi}
\frac{1}
{\theta\,{\phi^{1+\frac{2}{n}}}}
 \,\frac{\sin{\theta}}{\theta}
d \theta,
  \\
  \beta_{n} &:=&
\int_{0}^{\pi}
\frac{\,\phi^{1-\frac{2}{n}}(\theta)\,\sin{\theta}}
{\left(\pi-\theta\right)^{2}
\,\phi^{2}\left(\pi-\theta\right)
}\,d \theta,
  \\
  \gamma_{n} &:=&
A\int_{0}^{\pi} \phi^{1-\frac{2}{n}}\,\sin{\theta}
d \theta.
\end{eqnarray*}
Observe that $\lim_{n\rightarrow +\infty}
a_{n}=1.\,$ We shall show that, while $\lim_{n\rightarrow +\infty}
\alpha_{n}=+\infty,\,$ the sequences $\left\{\beta_{n}\right\}_{n\geq 1}$ and
$\left\{\gamma_{n}\right\}_{n\geq 1}$ are both convergent. Using this in (\ref{lim})
proves (\ref{limq}). \\
\indent Exploiting the continuity and positivity of
$\,{\sin{\theta}}/
\left({\theta^2\,{\phi^{1+\frac{2}{n}}}}\right)\,$
on $[\pi/2,\pi]$, then applying the inequality
$\,{\sin{\theta}}/{\theta}\geq {2}/{\pi}\,$
when $\,0\leq \theta \leq \pi/2,\,$ we obtain
\begin{eqnarray}
\nonumber \alpha_{n}&=&
\int_{0}^{\pi/2}
\frac{1}
{\theta\,{\phi^{1+\frac{2}{n}}}}
 \,\frac{\sin{\theta}}{\theta}
d \theta+
\int_{\pi/2}^{\pi}
\frac{\sin{\theta}}
{\theta^2\,{\phi^{1+\frac{2}{n}}}}
 \,d \theta\\
\label{alpha11}&\geq&\frac{2}{\pi}
\int_{0}^{\pi/2}
\frac{1}
{\theta\,{\phi^{1+\frac{2}{n}}}}
 \,d \theta\,=\,\frac{n}{\pi (1+\log(2))^{\frac{2}{n}}}.
 \end{eqnarray}
This proves the divergence of $\{\alpha_{n}\}$.
Next, by the dominated convergence theorem
and (\ref{calc1}) we readily find
\begin{align*}
\lim_{n\rightarrow +\infty}  \gamma_{n} \,=\,
A\lim_{n\rightarrow +\infty} \int_{0}^{\pi} {\phi^{1-\frac{2}{n}}(\theta)} \,\sin{\theta}
d\theta\,=\,\int_{0}^{\pi} {\phi(\theta)} \,\sin{\theta} d\theta\,\lesssim 1.
\end{align*}
Finally, since $\displaystyle \,\theta \mapsto{\sin{\theta}}/{\left(\left(\pi-\theta\right)^{2}
\,\phi^{2}\left(\pi-\theta\right)\right)
}\in C\left([0,\pi/2]\right),\,$ then
using the local integrability of $\phi$ and
the dominated convergence theorem again implies
\begin{equation}\label{beta1}
\lim_{n\rightarrow \infty}
\int_{0}^{\pi/2}
\frac{\phi^{1-\frac{1}{n}}(\theta)\,\sin{\theta}}
{\left(\pi-\theta\right)^{2}
\,\phi^{2}\left(\pi-\theta\right)
}\,d \theta\,=\,\int_{0}^{\pi/2}
\frac{\phi(\theta)\,\sin{\theta}}
{\left(\pi-\theta\right)^{2}
\,\phi^{2}\left(\pi-\theta\right)
}\,d \theta\,\lesssim\,1.
\end{equation}
Furthermore, since $\,\phi\in C\left([\pi/2,\pi]\right),\,$ and
$\,\displaystyle \frac{\sin{\theta}}{\pi-\theta}=
\frac{\sin{\left(\pi-\theta\right)}}{\pi-\theta}\leq 1,\,$ on $[\pi/2,\pi],\,$ then
\begin{equation}\label{beta2}
\int_{\pi/2}^{\pi}
\frac{\phi^{1-\frac{1}{n}}(\theta)\,\sin{\theta}}
{\left(\pi-\theta\right)^{2}
\,\phi^{2}\left(\pi-\theta\right)
}\,d \theta  \,\lesssim\,
\int_{\pi/2}^{\pi}
\frac{d \theta}
{\left(\pi-\theta\right)
\,\phi^{2}\left(\pi-\theta\right)}\,\approx\,1.
\end{equation}
The convergence of $\{\beta_{n}\}$
follows from (\ref{beta1}) together with
(\ref{beta2}). \\\\
\indent The proof of (\ref{main22})
shares the main idea of (\ref{main21}). The functions
$\,g_{n}(\theta,\varphi):=
{\psi(\theta)}^{\frac{1}{2}-\frac{1}{n}}\in
L^{2}\left(\mathbb{S}^{2}\right),\,n\geq 1,\,$
and satisfy $\,\displaystyle \lim_{n\rightarrow \infty}\frac{Q\left(g_{n};\psi\right)}
{T_{B}(g_{n})}=1$. Indeed, we have
\begin{eqnarray*}
\frac{Q\left(g_{n};\psi\right)}
{T_{B}(g_{n})}&=&
\frac{
\int_{0}^{\pi}
\frac{g_{n}^{2}\,\sin{\theta}}{\theta^{2}\,\psi^{2}(\theta) }\,d\theta+
\int_{0}^{\pi}
\frac{g_{n}^{2}\,\sin{\theta}}
{\left(\pi-\theta\right)^{2}
\,\psi^{2}\left(\pi-\theta\right)
}\,d\theta}
{4\int_{0}^{\pi} \left(\frac{\partial g_{n}}{\partial \theta} \right)^{2}\,\sin{\theta}
d\theta+B \int_{0}^{\pi} g_{n}^{2} \,\sin{\theta}
d\theta}\\
&=&\frac{\tilde{\alpha}_{n}}
{a_{n}\tilde{\alpha}_{n}+\tilde{\beta}_{n}}=
\frac{1}{a_{n}}\left(
1-\frac{\tilde{\beta}_{n}/a_{n}}{\tilde{\alpha}_{n}+
\tilde{\beta}_{n}/a_{n}}\right)
\end{eqnarray*}
where
\begin{eqnarray*}
  \tilde{\alpha}_{n} &:=&\int_{0}^{\pi}
\frac{\sin{\theta}\,d \theta}
{\theta^2\,{\psi^{1+\frac{2}{n}}}}
 +\int_{0}^{\pi}
\frac{\sin{\theta}d \theta}
{(\pi-\theta)^2\,{\psi^{1+\frac{2}{n}}}}
=2\int_{0}^{\pi}
\frac{\sin{\theta}\,d \theta}
{\theta^2\,{\psi^{1+\frac{2}{n}}}},
  \\
\tilde{\beta}_{n} &:=&
B\int_{0}^{\pi} \psi^{1-\frac{2}{n}}\,\sin{\theta}
d \theta-
a_{n}\int_{0}^{\pi}
M(\theta)
\frac{\sin{\theta}}{\psi^{1+\frac{2}{n}}}
\,d \theta.
\end{eqnarray*}
Similarly to (\ref{alpha11}), we have
\begin{eqnarray*}
\tilde{\alpha}_{n}
&=&2\int_{0}^{1}
\frac{\sin{\theta}}{\theta^{2}}
\frac{1}{\,{\psi^{1+\frac{2}{n}}}}\,d \theta+
2\int_{1}^{\pi}
\frac{\sin{\theta}}{\theta^{2}}
\frac{1}{\,{\psi^{1+\frac{2}{n}}}}\,d \theta\\
&\geq&
2\int_{0}^{1}
\frac{\sin{\theta}}{\theta^{2}}
\frac{1}{\,{\psi^{1+\frac{2}{n}}}}\,d \theta
\,=\,2\int_{0}^{1}
\frac{\sin^{2}{\theta}}{\theta^{2}\cos{\theta}}
\frac{1}{\,{\psi^{1+\frac{2}{n}}}}
\frac{\cos{\theta}}{\sin{\theta}}
\,d \theta\\
&\geq&\frac{8}{\pi^{2}}
\int_{0}^{1}
\frac{1}{\,{\psi^{1+\frac{2}{n}}}}
\frac{\cos{\theta}}{\sin{\theta}}\,d \theta
\,=\,\frac{4n}{\pi^2}
\frac{1}{\left(1+\log{\pi}\right)^{\frac{2}{n}}}.
\end{eqnarray*}
Hence $\,\lim_{n\rightarrow \infty} \tilde{\alpha}_{n}=
\infty$. Recall from (\ref{mm1}) and
(\ref{mm2}) that $M\in C([0,\pi])$. Also, since $\psi\in L^{1}_{\text{loc}} \left(\mathbb{R}\right),\,$ $\psi> 1$ uniformly, then $\,\lim_{n\rightarrow \infty} \tilde{\beta}_{n}$ exists by the dominated convergence theorem.
\end{proof}
\section{Proof of Theorem \ref{thm3}}
\begin{proof}
The transition to the inequalities (\ref{main11}) and (\ref{main12}) from their respective
stronger versions, (\ref{p4}) and (\ref{p5}), comes from the bounds
 \begin{equation*}
\int_{\mathbb{S}^{2}} F(\theta)\,  f^{2} d\sigma\,\leq\, A\int_{\mathbb{S}^{2}} f^{2} d\sigma,\quad
\int_{\mathbb{S}^{2}} G(\theta)\,  f^{2} d\sigma\,\leq\, B\int_{\mathbb{S}^{2}} f^{2} d\sigma
\end{equation*}
where the bounded functions $F$ and $G$ are both positive and independent of $f$.
Interestingly, as seen in Section \ref{proofthm2}, the size of $\,0<A,B<\infty\,$ played no role
in optimising (\ref{main11}) and (\ref{main12}).\\
\indent Up to the inequality (\ref{p4}) or (\ref{p5})
an equality relation persists except for the only inequality (\ref{notattain}). So
a sufficient and necessary condition for an equality
in (\ref{p4}) or (\ref{p5}) (and a necessary condition for an equality in (\ref{main11}) and (\ref{main12}))
is an equality in (\ref{notattain}). But an equality
in (\ref{notattain}) occurs if and only if
\begin{equation}\label{iff1}
\int_{\mathbb{S}^{2}}
 \psi |\nabla  \,  g|^{2}
  d\sigma=0.
\end{equation}
Recalling that $g=f/\sqrt{\psi},$ we compute
\begin{eqnarray}
\nonumber
\psi|\nabla   g|^{2}&=&
\psi\left|\frac{\nabla f}{\sqrt{\psi}}-
\frac{1}{2}\frac{f}{\psi^{\frac{3}{2}}}\frac{\partial \psi}{\partial \theta} \hat{\theta}
\right|^{2}
\\
\nonumber  &=&|\nabla f|^{2}-
\frac{f}{\psi}\frac{\partial \psi}{\partial \theta} \, \nabla f\cdot\hat{\theta}+
\frac{1}{4}\frac{f^2}{\psi^{2}}\left(\frac{\partial \psi}{\partial \theta}\right)^{2}\\
\nonumber  &=&|\nabla f|^{2}-\left(\frac{\partial f}{\partial \theta}\right)^{2}
+\left(\frac{\partial f}{\partial \theta}\right)^{2}
-\frac{f}{\psi}\frac{\partial \psi}{\partial \theta} \, \frac{\partial f}{\partial \theta}+
\frac{1}{4}\frac{f^2}{\psi^{2}}\left(\frac{\partial \psi}{\partial \theta}\right)^{2}\\
\label{iff2}&=&|\nabla f|^{2}-\left(\frac{\partial f}{\partial \theta}\right)^{2}
+\left(\frac{\partial f}{\partial \theta}
-\frac{1}{2}\frac{f}{\psi}\frac{\partial \psi}{\partial \theta}\right)^{2}.
\end{eqnarray}
Since $\,\displaystyle |\nabla f|^{2}-\left(\frac{\partial f}{\partial \theta}\right)^{2}=\frac{1}{\sin^{2}{\theta}}
\left(\frac{\partial f}{\partial \varphi}\right)^{2}\,\geq\,0,\,$ then,
by (\ref{iff2}), the equality (\ref{iff1})
is equivalent to
\begin{equation}\label{iff3}
\int_{\mathbb{S}^{2}}
|\nabla f|^{2}-\left(\frac{\partial f}{\partial \theta}\right)^{2}d\sigma\,=\,
\int_{\mathbb{S}^{2}}
\left(\frac{\partial f}{\partial \theta}
-\frac{1}{2}\frac{f}{\psi}\frac{\partial \psi}{\partial \theta}\right)^{2}d\sigma\,=\,0.
\end{equation}
The equalities (\ref{iff3}) are, in their turn,
equivalent to
\begin{equation}\label{iff4}
\frac{1}{\sin{\theta}}
\left|\frac{\partial f}{\partial \varphi}\right|\,=\,
\left|\frac{\partial f}{\partial \theta}
-\frac{1}{2}\frac{f}{\psi}\frac{\partial \psi}{\partial \theta}\right|\,=\,0.
\end{equation}
Suppose that $f$ is not the zero function. Then (\ref{iff4}) are possible if and only if
\begin{equation*}
f\,=\,f(\theta),\,
\frac{d f}{f}
\,=\,\frac{1}{2}\frac{d\psi}{\psi}.
\end{equation*}
That is $f\,=\,c\sqrt{{\psi}},$\, $c$ is a constant.
But such $f\notin H^{1}\left(\mathbb{S}^{2}\right)\,$ because
\begin{eqnarray*}
\int_{\mathbb{S}^{2}}
|\nabla  \,  f|^{2}
  d\sigma&=&2\pi\int_{0}^{\pi}
\left(\frac{\partial f}{\partial \theta}\right)^{2}
  d\theta\,\gtrsim
\int_{0}^{1}
\frac{\cos^{2}{\theta}}{\sin{\theta}}
\frac{1}{\psi}  d\theta
\\
&\gtrsim&\int_{0}^{1}
\frac{d\theta}{\sin{\theta}\,\phi(\sin{\theta})}
\,\approx\,
\int_{0}^{1}
\frac{d\theta}{{\theta}\,\phi({\theta})}
\,=\,+\infty.
\end{eqnarray*}
\end{proof}
\section{Proof of Theorem \ref{thm4}}
\begin{proof}
Write
\begin{align*}
\frac{1}{\theta}\frac{1}{\phi^{2}\left(\theta\right)} =
\nabla \left( \frac{1}{\phi\left(\theta\right)}\right)\cdot \hat{\theta}.
\end{align*}
Assume that $f$ is smooth. Then integrating by parts w.r.t. the surface measure $\sigma$ we get
\begin{align}
\nonumber  \int_{\mathbb{S}^{2}} \frac{f^{2}}{\theta^{2}\phi^{2}\left(\theta\right)}
\, d\sigma =&  \int_{\mathbb{S}^{2}} \nabla \left( \frac{1}{\phi\left(\theta\right)}\right)\cdot \frac{f^{2}}{\theta} \hat{\theta}  d\sigma \\
\nonumber   =&- \int_{\mathbb{S}^{2}}  \frac{1}{\phi\left(\theta\right)} \nabla \cdot \left(\frac{f^{2}}{\theta} \hat{\theta}\right)d\sigma\\
\nonumber   =&- 2\int_{\mathbb{S}^{2}}  \frac{f\,\nabla f.\hat{\theta}}{\theta\,\phi\left(\theta\right)}d\sigma+
\int_{\mathbb{S}^{2}}  \frac{f^{2}}{\theta^{2}\,\phi\left(\theta\right)}d\sigma+\\
\label{w1} &-
\int_{\mathbb{S}^{2}}  \frac{f^{2}}{\theta\,\phi\left(\theta\right)}
\frac{\cos{\theta}}{\sin{\theta}}d\sigma.
\end{align}
Observe here that each of the last two integrals on the right hand side of (\ref{w1}) can diverge. They suffer  nonintegrable singularities at $\theta=0$.
The reality is, put together, their sum
\begin{align}\label{wc1}
I:=\int_{\mathbb{S}^{2}}  \frac{f^{2}}{\theta^{2}\,\phi\left(\theta\right)}d\sigma-
\int_{\mathbb{S}^{2}}  \frac{f^{2}}{\theta\,\phi\left(\theta\right)}
\frac{\cos{\theta}}{\sin{\theta}}d\sigma=
\int_{\mathbb{S}^{2}} \frac{1}{\theta\,\phi\left(\theta\right)}
\left(\frac{1}{\theta}-\frac{\cos{\theta}}{\sin{\theta}} \right)f^{2}d\sigma
\end{align}
is convergent. In fact
\begin{align*}
\lim_{\theta\rightarrow 0^{+}}\frac{1}{\theta\,\phi\left(\theta\right)}
\left(\frac{1}{\theta}-\frac{\cos{\theta}}{\sin{\theta}} \right)=0.
\end{align*}
Also,$\,\theta\mapsto1/\left({\theta^{2}\,\phi\left(\theta\right)}\right)\,$
is continuous on a neighborhood of $\theta=\pi$.
Furthermore, if we fix $\delta>0$ and let
$\,D:=\left\{x(\theta,\varphi)\in \mathbb{S}^{2}:
0\leq \theta < \delta \right\}$,
then the integral $
\displaystyle \int_{\mathbb{S}^{2}\setminus D} \frac{f^{2}}{\theta\,\phi\left(\theta\right)}
\frac{\cos{\theta}}{\sin{\theta}} d\sigma$
does exist.
Unfortunately, we can not control the integral $I$ by $\int_{\mathbb{S}^{2}} f^{2} d\sigma$, up to a constant factor. The reason is
\begin{align*}
\lim_{\theta\rightarrow \pi^{-}}\frac{1}{\theta\,\phi\left(\theta\right)}
\frac{\cos{\theta}}{\sin{\theta}} =\infty.
\end{align*}
But since
\begin{align*}
\lim_{\theta\rightarrow \pi^{-}}\left(\frac{1}{\theta\,\phi\left(\theta\right)}
\frac{\cos{\theta}}{\sin{\theta}}+
\frac{1}{\pi}\frac{1}{(\pi-\theta)}\right)
 =0
\end{align*}
then, we may introduce the convergent integral  $\displaystyle J:=\frac{1}{\pi}\int_{\mathbb{S}^{2}}
\frac{f^{2}}{\pi-\theta}\, d\sigma$
to the integral $I$ to get
\begin{equation}\label{wc2}
I=I-J+J=\int_{\mathbb{S}^{2}}
K(\theta)\,f^{2}\, d\sigma+J
\end{equation}
where
\begin{align*}
K(\theta):=\frac{1}{\theta\,\phi\left(\theta\right)}
\left(\frac{1}{\theta}-\frac{\cos{\theta}}{\sin{\theta}} \right)-\frac{1}{\pi}\frac{1}{(\pi-\theta)}.
\end{align*}
By the continuity of $K$ on $]0,\pi[$ and since
\begin{align*}
\lim_{\theta\rightarrow 0^{+}}
K(\theta)=-\lim_{\theta\rightarrow \pi^{-}}
K(\theta)=-\frac{1}{\pi^{2}}
\end{align*}
then $K$ is bounded on $[0,\pi]$. Actually,
$K$ is monotonically increasing. Thus
\begin{align}\label{wc3}
\sup_{[0,\pi]}|K|=\frac{1}{\pi^{2}}.
\end{align}
Using (\ref{wc3}) in (\ref{wc2}) we deduce that
\begin{align}\label{wc4}
I\leq \frac{1}{\pi^{2}}\int_{\mathbb{S}^{2}}
f^{2}\, d\sigma+J.
\end{align}
Returning with (\ref{wc4}) to the inequality (\ref{w1})
in the light of (\ref{wc1}) we obtain
\begin{equation}
\label{w2}  \int_{\mathbb{S}^{2}} \frac{f^{2}}{\theta^{2}\phi^{2}\left(\theta\right)}
\, d\sigma\,\leq\,- 2\int_{\mathbb{S}^{2}}  \frac{f\,\nabla f.\hat{\theta}}{\theta\,\phi\left(\theta\right)}d\sigma+
\frac{1}{\pi^{2}}\int_{\mathbb{S}^{2}}
f^{2}\, d\sigma+\frac{1}{\pi}\int_{\mathbb{S}^{2}}
\frac{f^{2}}{\pi-\theta}\, d\sigma.
\end{equation}
Applying Cauchy's inequality with an $\epsilon$
we find
\begin{equation}\label{w4}
- 2\int_{\mathbb{S}^{2}}  \frac{f\,\nabla f.\hat{\theta}}{\theta\,\phi\left(\theta\right)}d\sigma
\leq 2\epsilon \int_{\mathbb{S}^{2}}
\frac{f^{2}}{\theta^{2}\phi^{2}\left(\theta\right)}
\, d\sigma+\frac{1}{2\epsilon}\int_{\mathbb{S}^{2}}
|\hat{\theta}.\nabla f|^{2} d\sigma.
\end{equation}
Therefore, it follows from (\ref{w2})
and (\ref{w4}) that
\begin{align}
\nonumber
&2\epsilon(1-2\epsilon)
\int_{\mathbb{S}^{2}} \frac{f^{2}}{\theta^{2}\phi^{2}\left(\theta\right)}
\, d\sigma \leq \int_{\mathbb{S}^{2}}
|\hat{\theta}.\nabla f|^{2} d\sigma+\\
\label{w5}&\qquad\qquad\qquad
+\frac{2\epsilon}{\pi^{2}}\int_{\mathbb{S}^{2}}
f^{2}\, d\sigma+\frac{2\epsilon}{\pi}\int_{\mathbb{S}^{2}}
\frac{f^{2}}{\pi-\theta}\, d\sigma,\quad 0<\epsilon<\frac{1}{2}.
\end{align}
The choice $\epsilon={1}/{4}$ maximizes
the factor $2\epsilon(1-2\epsilon)$ and,
consequently, the left hand side of (\ref{w5}).
This proves (\ref{w01}).
The inequality (\ref{w02}) can be obtained analogously.\\
\indent In the fashion of the proof
of Theorem \ref{thm2}, the sequence
$f_{n}=\phi^{\frac{1}{2}-\frac{1}{n}}$
clearly satisfies
\begin{equation*}
\lim_{n\rightarrow \infty}
\frac{\frac{1}{4}\int_{0}^{\pi} \frac{f_{n}^{2}}{{{\rho}_{\phi}}^{2}
\left(\theta\right)}
\, \sin{\theta}\,d \theta}
{U(f_{n}
)+\frac{1}{2\pi}\int_{0}^{\pi}
\frac{f_{n}^{2}}{\pi-\theta}\, \sin{\theta}\,d \theta} =
\lim_{n\rightarrow \infty}
\frac{\frac{1}{4}\int_{0}^{\pi}
\frac{f_{n}^{2}}{{{\rho}_{\phi}}^{2}
\left(\pi-\theta\right)}
\, \sin{\theta}\,d \theta}
{U(f_{n})+\frac{1}{2\pi}\int_{0}^{\pi}
\frac{f_{n}^{2}}{\theta}\, \sin{\theta}\,d \theta} =
1
\end{equation*}
where
\begin{equation*}
 U(f)=\int_{0}^{\pi}
\left(\frac{\partial f}{\partial{\theta}}\right)^{2}\, \sin{\theta} \,d \theta+\frac{1}{2\pi^{2}}\int_{0}^{\pi}
f^{2}\, \sin{\theta}\,d \theta.
\end{equation*}
One only needs to inspect the convergence of
$\,\int_{0}^{\pi}
\left(
\phi^{1-\frac{2}{n}}\sin{\theta}/{\theta}\,\right)d \theta,$\\
$\;\int_{0}^{\pi}
\left(\phi^{1-\frac{2}{n}}\sin{\theta}/\left(\pi-{\theta}\right)
\right)d \theta
$ as $n\rightarrow \infty.\,$
This is obvious from the bound
$\sin{\theta}\leq \min\{{\theta},{\pi-\theta}\}\,$
on $[0,\pi]$ and the fact $\phi \in L^{1}\left([0,\pi]\right)$.\\
\indent Finally, careful review of the proof of (\ref{w01})
above reveals that a necessary condition
for a function $f\in H^{1}\left(\mathbb{S}^{2}\right)\setminus\{0\}$
to achieve an equality in (\ref{w01}) is
that it yields an equality in (\ref{w4}).
This is equivalent to
\begin{equation}\label{ww}
\nabla f.\hat{\theta}\,=\,
- \frac{1}{2}\frac{f}{\theta\,\phi\left(\theta\right)}.
\end{equation}
Suppose (\ref{ww}) was true. Then
by (\ref{w1}) and (\ref{wc1})
we must have
\begin{equation}\label{ww1}
\int_{\mathbb{S}^{2}} \frac{h(\theta)\,f^{2}}{\theta\,\phi\left(\theta\right)}
\,d\sigma\,=\,0
\end{equation}
where
\begin{equation*}
 h(\theta):=\frac{1}{\theta}-\frac{\cos{\theta}}{\sin{\theta}}.
\end{equation*}
On the other hand
\begin{equation*}
\lim_{\theta\rightarrow 0^{+}} h(\theta)=0,\;\;
 h^{\prime}(\theta)\,=\,
\frac{{\theta}^{2}-
\sin^{2}{\theta}}{{\theta}^{2}\,\sin^{2}{\theta}}>0,\;
0<\theta<\pi.
\end{equation*}
This shows $h$ is strictly positive on $]0,\pi]$ and since $\,\theta \phi(\theta)\geq 0$  then (\ref{ww1}) is a contradiction.
\end{proof}



\textbf{References}
\bibliographystyle{model1a-num-names}
\bibliography{<your-bib-database>}

\end{document}